\documentclass[10pt]{amsart}

\usepackage{graphicx}        
\usepackage{multicol}        
\usepackage[bottom]{footmisc}
\usepackage{amscd,amsmath,amssymb,amsfonts}
\usepackage[cmtip, all]{xy}

\usepackage{hyperref}

\usepackage{bbm}
\usepackage{tikz-cd}
\usepackage{enumerate}
\usepackage{mathrsfs}
\usepackage{amsthm}
\usepackage{comment}

\newtheorem{theorem}{Theorem}[section] 
\newtheorem{proposition}[theorem]{Proposition}
\newtheorem{lemma}[theorem]{Lemma}

\theoremstyle{definition}
\newtheorem{definition}[theorem]{Definition}
\theoremstyle{remark}

\numberwithin{equation}{section}

\newcommand{\CH}{{\operatorname{CH}}}

\newcommand{\Pic}{{\rm Pic}}
\newcommand{\sPic}{{\operatorname{\mathcal{P}ic}}}

\newcommand{\Hom}{{\rm Hom}}

\newcommand{\Spec}{{\rm Spec\,}}

\newcommand{\Tr}{{\text{Tr}}}

\newcommand{\sE}{{\mathcal E}}
\newcommand{\sF}{{\mathcal F}}
\newcommand{\sG}{{\mathcal G}}

\newcommand{\sI}{{\mathcal I}}

\newcommand{\sK}{{\mathcal K}}
\newcommand{\sL}{{\mathcal L}}
\newcommand{\sM}{{\mathcal M}}

\newcommand{\sO}{{\mathcal O}}
\newcommand{\sP}{{\mathcal P}}

\newcommand{\sU}{{\mathcal U}}
\newcommand{\sV}{{\mathcal V}}
\newcommand{\sW}{{\mathcal W}}

\newcommand{\A}{{\mathbb A}}

\newcommand{\C}{{\mathbb C}}

\renewcommand{\L}{{\mathbb L}}

\renewcommand{\P}{{\mathbb P}}

\newcommand{\Z}{{\mathbb Z}}

\newcommand{\BM}{{\operatorname{B.M.}}}

\renewcommand{\det}{\operatorname{det}}
\newcommand{\Nm}{{\operatorname{Nm}}}

\newcommand{\id}{{\operatorname{\rm Id}}}

\renewcommand{\dim}{{\operatorname{\rm dim}}}

\newcommand{\Hilb}{\operatorname{Hilb}}

\newcommand{\del}{\partial}

\newcommand{\Sym}{{\operatorname{Sym}}}

\newcommand{\Gr}{{\operatorname{\rm Gr}}} 
\newcommand{\rnk}{{\operatorname{rank}}}

\newcommand{\KQ}{{\operatorname{KQ}}}

\newcommand{\GW}{{\operatorname{GW}}} 
 
\newcommand{\SH}{{\operatorname{SH}}}

\newcommand{\sHom}{\mathcal{H}om}

\newcommand{\SL}{\operatorname{SL}}

\newcommand{\Ass}{\operatorname{Ass}}

\newcommand{\vir}{\text{\it vir}}
\newcommand{\perf}{\text{\it perf}}

\newcommand{\ind}[1]{}

\newcommand{\inp}[1]{}

\newcommand{\EM}{{\operatorname{EM}}}

\newcommand{\Cone}{{\operatorname{Cone}}}

\newcommand{\RsHom}{{{R\mathcal H}om}}

\makeatletter
 \@namedef{subjclassname@2020}{\textup{2020} Mathematics Subject Classification}
\makeatother

\date{ \today}

\author[M.~Levine]{Marc~Levine}
\address{Marc~Levine, Universit\"at Duisburg-Essen,
Fakult\"at Mathematik, Campus Essen, 45117 Essen, Germany}
\email{marc.levine@uni-due.de}

\thanks{This paper is part of a project that has received funding from the European Research Council (ERC) under the European Union's Horizon 2020 research and innovation programme (grant agreement No. 832833).\\ 
\includegraphics[scale=0.08]{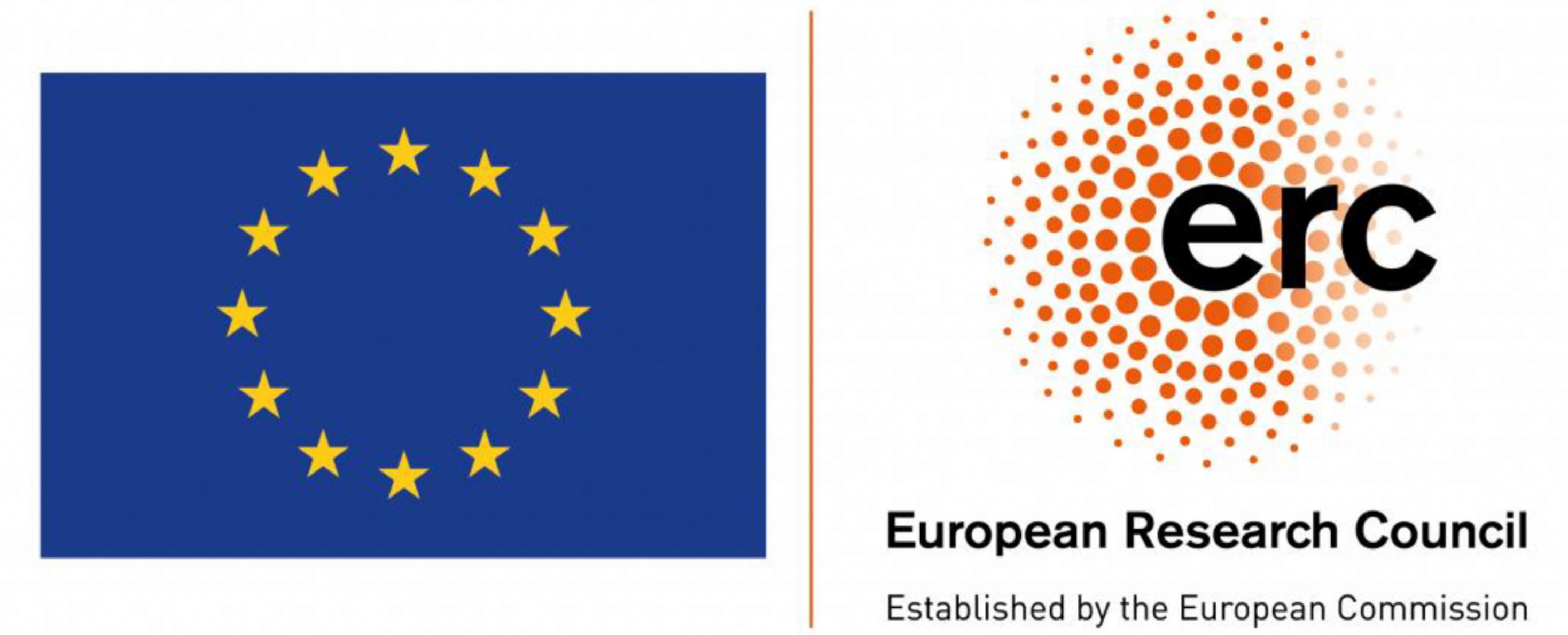}}

\title{Orienting the Hilbert scheme of points on a spin threefold}

\subjclass[2020]{
Primary 14N10, 14N35
Secondary 14Q30, 	14F42}
\keywords{enumerative geometry, quadratic refinements, Hilbert scheme, orientations}

\begin{document}

\begin{abstract}  We give an algebraic construction of orientation data for the Hilbert scheme of 0-dimensional subschemes of a smooth projective threefold endowed with a choice of square root of the canonical sheaf. 
\end{abstract}

\maketitle

\tableofcontents

\section{Introduction} In \cite{JU} Joyce and Upmeier construct orientations on certain moduli stacks for a Calabi-Yau threefold $X$ over $\C$, or more generally, for a smooth proper threefold $X$ over $\C$ endowed with a square root of the canonical sheaf: a {\em spin structure} on $X$. The orientation data for an associated moduli stack $\sM(X$) is similarly a isomorphism $K_\sM\cong \sL_\sM^{\otimes 2}$ for a chosen invertible sheaf $\sL_\sM$ on $\sM(X)$ (together with some additional compatibilities that we omit mentioning); here $K_\sM$ is the {\em virtual canonical sheaf} on the moduli stack. As mentioned in the introduction to {\it loc. cit.}, the existence of an orientation on $\sM(X)$ is needed for the various refinements of Donaldson-Thomas invariants, for example, a refinement to $K_0(Var_\C)$, constructed by Kontsevich-Soibelman \cite[\S 5]{KS}.  The construction used in \cite{JU} relies on both analytic and topological inputs, and does not appear to have an easy extension to the purely algebraic setting of a smooth projective threefold $X$ over a field $k$, endowed  a spin structure (also over $k$).

In \cite{Toda}, Y. Toda constructs an orientation on the Hilbert scheme of 0-dimensional subschemes of a Calabi-Yau threefold over $\C$, by a different method from that used in \cite{JU}. 

My interest in this question arose from my construction in \cite{LevineVFC} of a ``motivic'' virtual fundamental class $[Z,\phi]^\vir_\sE$ associated to a perfect obstruction theory $\phi:E_\bullet\to \L_{Z/k}$ on a finite type $k$-scheme $X$, where $\L_{Z/k}$ is the truncation $\tau^{\ge-1}L_{Z/k}$ of the cotangent complex, and $\sE$ is a motivic ring spectrum, that is, $\sE$ is an object in the motivic stable homotopy category over $k$, $\sE\in \SH(k)$, endowed with a commutative, associative, unital multiplication $\sE\wedge\sE\to \sE$.  Assuming $\sE$ has a so-called $\SL$-orientation \cite{Ananyevskiy, Panin}, the virtual fundamental class lives in the twisted motivic Borel-Moore homology:
\[
[Z,\phi]^\vir_\sE\in \sE_{2r,r}^\BM(Z,\det E_\bullet)
\]
where $r$ is virtual rank of $E_\bullet$ and $\det E_\bullet$ is the virtual determinant. If $Z$ is proper over $k$, $r=0$ and we give an isomorphism $\rho:\det E_\bullet\xrightarrow{\sim} \sM^{\otimes 2}$ for some invertible sheaf $\sM$ on $Z$, we have a degree map
\[
\deg^{\sE, \rho}_k:\sE_{2r,r}^\BM(Z,\det E_\bullet)\to \sE^{0,0}(\Spec k).
\]

In our envisaged application, $Z=\Hilb^n(X)$, the Hilbert scheme of length $n$ 0-dimensional closed subschemes of a smooth projective threefold $X$, endowed with a spin structure $\rho_X:\omega_{X/k}\xrightarrow{\sim}\sL^{\otimes2}$. One has the rank zero perfect obstruction theory $\phi_T:E_\bullet\to  L_{\Hilb^n(X)/k}$ given by Thomas \cite{Thomas} and Behrend-Fantechi \cite[Lemma 1.2.2]{BehrendFantechi}, and our main result here is to give an invertible sheaf $\sM$ on $\Hilb^n(X)$ and canonical isomorphism $\rho_{\Hilb^n}:\det E_\bullet\xrightarrow{\sim} \sM^{\otimes 2}$ (depending on the choice of spin structure $\rho_X:\omega_{X/k}\xrightarrow{\sim}\sL^{\otimes2}$). The canonical class $\sK_\sM$ from Joyce-Upmeier is just $\det E_\bullet$ in this case, so we refer to an isomorphism of the form 
$\rho:\det E_\bullet\xrightarrow{\sim} \sM^{\otimes 2}$ as a {\em virtual orientation} of the obstruction theory $E_\bullet$. In case the perfect obstruction theory arises from a quasi-smooth derived scheme $\sM$,  $\det E_\bullet$ is the canonical sheaf on $\sM$, so the notation is more or less consistent.

We then take $\sE$ to be a motivic spectrum such that $\sE^{0,0}(k)$ is isomorphic to 
the Grothendieck-Witt ring $\GW(k)$ of virtual non-degenerate symmetric bilinear  forms over $k$. In such a case,  our virtual orientation $\rho_{\Hilb^n}$ gives a quadratic DT invariant
\[
\deg_k^\sE([\Hilb^n(X), \phi_T]^\vir_\sE)\in \GW(k).
\]
Typical choices for $\sE$ are  include Hermitian $K$-theory $\KQ$, Chow-Witt motivic cohomology $\widetilde{H\Z}$ or the Eilenberg-MacLane spectrum $\EM(\sK^{MW}_*)$ associated to Milnor-Witt $K$-theory $\sK^{MW}_*$ (see
Morel's theorem \cite[Theorem 6.2, Theorem 6.2.2]{MorelPi0}).

In all these cases, the virtual fundamental class and the degree live in the ``geometric part''  of $\sE$. For Chow-Witt motivic cohomology, these are the Chow-Witt groups $\widetilde{\CH}_*(-)$. One can also use a somewhat simpler theory, $H^*(-, \sW)$, the cohomology of the sheaf of Witt rings and the resulting degree map lands in $H^0(\Spec k, \sW)=W(k)$, the Witt ring of $k$. One has maps $\widetilde{\CH}^*(-)\to \CH^*(-)$ (the usual Chow groups) and $\widetilde{\CH}^*(-)\to H^*(-\sW)$, 
transforming the virtual fundamental class $[\Hilb^n(X), \phi_T]^\vir_{\widetilde{\CH}}$ to the classical one in the Chow groups $[\Hilb^n(X), \phi_T]^\vir_\CH$ and to the Witt-sheaf version 
$[\Hilb^n(X), \phi_T]^\vir_{H^*(-,\sW)}$, respectively. One has a similar compatibility for the degrees, i.e., the corresponding DT invariants, so 
\[
\deg_k^{\widetilde{\CH}}([\Hilb^n(X), \phi_T]^\vir_{\widetilde{\CH}})\in \GW(k)
\]
maps to the usual DT invariant $\deg_k^{{\CH}}([\Hilb^n(X), \phi_T]^\vir_{{\CH}})\in\Z$ under the rank map $\GW(k)\to \Z$, and to the Witt-sheaf variant $\deg_k^{\sW}([\Hilb^n(X), \phi_T]^\vir_\sW\in W(k)$ under the canonical surjection $\GW(k)\to W(k)$. 

This paper is organized as follows. We first discuss some commutative algebra that gives a sufficient condition for a restriction map $j^*:\Pic(X)\to \Pic(U)$ associated to an open immersion $j$ to be injective in \S\ref{sec:Triv}. In \S\ref{sec:Det} we briefly recall the construction of the norm map on $\Pic$ and a result of Deligne relating this to the determinant of the pushforward of a perfect complex under a finite flat map. In \S\ref{sec:Orient} we give our modification of Toda's argument to compute $\det E_\bullet$ for the Thomas perfect obstruction theory on $\Hilb^n$, and we state and prove our main theorem in \S\ref{sec:Main}

For computations, it is useful to introduce a group action for a suitable group and use the form of Atiyah-Bott localization and the virtual Bott residue theorem as constructed in \cite{LevineVirBott}. For this, we briefly describe how to extend our construction of the orientation for $\Hilb^n$ to the equivariant setting in the last section of this paper, \S\ref{sec:equi}.

Yokinubo Toda kindly explained the proof of his result \cite[Proposition 3.1]{Toda} to me in response to my questions about his paper. His explanations  were very useful in developing this note and I wish to express my gratitude to him for these discussions. I also want to thank Andrea Ricolfo for explaining the results in his paper \cite{Ricolfi} on equivariant perfect obstruction theories,  and the use of the Atiyah class in constructing a perfect obstruction theory, as laid out in the paper of Huybrechts-Thomas \cite{HT}.

\section{Extending a trivialization of an invertible sheaf}\label{sec:Triv}

\begin{lemma}\label{lem1} Let $A\to B$ be a  homomorphism of noetherian rings, and take $f$ in $B$. Suppose that $B$ is smooth over $A$ and that, for each prime ideal $\mathfrak{p}$ of $A$, $f$ is non-zero on each irreducible component of $\Spec B/\mathfrak{p}B$. Then the map $B\to B[1/f]$ is injective.
\end{lemma}

\begin{proof} For a commutative ring $C$, we let $\Ass(C)$ denote the set of associated primes of $C$, that is the associated primes of the zero ideal in $C$. The map $B\to B[1/f]$ is injective if and only if $f$ is a non-zero divisor on $B$, equivalently, $f$ is not in any associated prime of $B$. Since $A\to B$ is flat we have by \cite[Chap. 3, Theorem 12]{MatsumuraComAlg}
\[
\Ass(B)=\cup_{\mathfrak{p}\in \Ass(A)}\Ass(B/\mathfrak{p}B).
\]
Since $B$ is smooth over $A$, the extension 
\[
A/\mathfrak{p}\to B\otimes_AA/\mathfrak{p}=B/\mathfrak{p}B
\]
is smooth over the integral domain $A/\mathfrak{p}$, so $\Spec B/\mathfrak{p}B$ is a disjoint union of its integral components. Equivalently, $B/\mathfrak{p}B$ is a product of integral domains
\[
B/\mathfrak{p}B=\prod_{i=1}^r\bar{B}_i
\]
By assumption, the image of $f$ in $\bar{B}_i$ is non-zero for each $i$, so $f$ is a non-zero divisor on $B/\mathfrak{p}B$, hence $f\not\in\mathfrak{q}$ for each $\mathfrak{q}\in \Ass(B)$, and thus $f$ is a non-zero divisor on $B$.
\end{proof}

\begin{lemma}\label{lem2} Let $(\sO, \mathfrak{m})$ be a Noetherian local ring, let $f,g$ be a regular sequence in $\mathfrak{m}$ and let $W\subset \Spec\sO$ be the closed subscheme with ideal $\sI_W=(f,g)\subset \sO$. Let $j_W:\Spec\sO\setminus W\to \Spec\sO$ be the inclusion, let $z\in \Spec\sO$ be the closed point,  and let $j_{W}^*:\sO\to (j_{W*}j^*_W\sO)_z$ be the canonical map. Then $j_{W}^*$ is an isomorphism.
\end{lemma}

\begin{proof} We cover $\Spec\sO\setminus W$ by $\Spec \sO[1/f]$ and $\Spec \sO[1/g]$. Taking $h\in   (j_{W*}j^*_W\sO)_z$, we have the restrictions $i_f^*h\in  \sO[1/f]$ and $i_g^*h\in  \sO[1/g]$, which we may write as $i_f^*h=h_1/f^n$, $i_g^*h=h_2/g^m$ for suitable $h_1, h_1\in \sO$, and integers $n,m\ge1$. Since $f,g$ is a regular sequence in $\mathfrak{m}$ by assumption, $g,f$ is also a regular sequence (this follows from \cite[Theorem 16.5]{MatsumuraComRingThy}), so both $f$ and $g$ are non-zero divisors in $\sO$. Thus the identity $h_1/f^n=h_2/g^m$ in $\sO[1/fg]$ implies that $h_1g^m=h_2f^n$ in $\sO$. 

On the other hand, $(f^n, g^m)$ is a regular sequence in $\mathfrak{m}$ by \cite[Theorem 16.1]{MatsumuraComRingThy}, so by \cite[Theorem 16.5]{MatsumuraComRingThy}, the Koszul complex
\[
0\to \sO\xrightarrow{(f^n, g^m)}\sO^2\xrightarrow{g^mp_1-f^np_2}\sO\to \sO/(f^n, g^m)\to 0
\]
is exact. Thus, the identity $h_1g^m=h_2f^n$ implies there is a unique $H\in \sO$ with $f^nH=h_1$, $g^mH=h_2$, which implies that $H/1=i_f^*h$ in $\sO[1/f]$ and $H/1=i_g^*h$ in 
$\sO[1/f]$. Since $\Spec \sO[1/f]\cup \Spec \sO[1/g]=\Spec\sO\setminus W$, this implies that $j_{W}^*H=h$, and thus $j_{W}^*$ is surjective. The uniqueness of $H$ implies that $j_{W}^*$ is injective, hence $j_{W}^*$ is an isomorphism, as claimed. 
\end{proof} 

\begin{lemma} \label{lem3} Let $Y$ be a smooth projective variety over $k$, and let $U$ be a finite type $k$-scheme. Let $Z$ be a closed subscheme of $U\times_kY$, equi-dimensional of dimension $d$ over $U$. Let $j:U\times_kY\setminus Z\to U\times_kY$ be the inclusion. Suppose that $d\le \dim_kY-2$. Then the canonical map
\[
j^*:\sO_{U\times_kY}\to j_*j^*\sO_{U\times_kY\setminus Z}
\]
is an isomorphism.
\end{lemma}

\begin{proof}  Clearly $j^*$ is an isomorphism on  $U\times_kY\setminus Z$, so we need only check on the stalks at a point $z\in Z$. Let $\sI_Z$ be the ideal sheaf of $Z$.  

Take $z\in Z$. Then 
there is an $f\in \sI_{Z,z}$ such that $f$ is non-zero on the irreducible component of $p_1(z)\times_k Y$ that contains $z$. For this, just take an affine open neighborhood $V=\Spec B$ of $z$ such that $p_1(z)\times Y\cap V$ is irreducible, take a point $y\in (p_1(z)\times Y\setminus Z)\cap V$ (which is non-empty by our condition on the dimension of $Z$ over $U$) and let $f$ be in $I_{Z\cap B}$ with $f(y)\neq 0$. By  Lemma~\ref{lem1}, the map $\sO_{U\times_kY, z}\to \sO_{U\times_kY, z}[1/f]$ is injective, as this factors through $j_z^*:\sO_{U\times_kY, z}\to (j_*j^*\sO_{U\times_kY})_z$, we see that $j_z^*$ is injective for all $z\in Z$, hence $j^*$ is injective.

Next, we show that $j^*_z$ is surjective for all $z\in Z$. For this, it suffices by Lemma~\ref{lem2} to find a regular sequence $f,g$ in $\sI_{Z,z}$.  Indeed, supposing we have done this, let $W=V(f,g)\subset \Spec\sO_{U\times_kY,z}$ and let $j_W:\Spec\sO_{U\times_kY, z}\setminus W\to \sO_{U\times_kY, z}$ be the inclusion. Then by Lemma~\ref{lem1}, the map
\[
j_W^*:\sO_{U\times_kY, z}\to (j_{W*}j_W^*\sO_{U\times_kY})_z
\]
is an isomorphism. Since $j_W^*$ factors as
\[
\sO_{U\times_kY, z}\xrightarrow{j^*} (j_*j^*\sO_{U\times_kY\setminus Z})_z\to  (j_{W*}j_W^*\sO_{U\times_kY})_z
\]
and $j^*$ is injective, it follows that $j^*$ is an isomorphism. 

Thus, it remains to see that we can find a regular sequence $f,g$ in $\sI_{Z,z}$.  We first do this in case   $U$ is smooth over $k$.  We show in fact that we can find a regular sequence $f,g$ in $\sI_{Z,z}$ such that the subscheme $W\subset \Spec \sO_{U\times_kY,z}$ with ideal $(f,g)$ is flat over $\sO_{U, p_1(z)}$. For this, we consider the global setting $Z\subset U\times_kY$; we may assume that $U=\Spec A$  with $A$ a finitely generated $k$-algebra. Let $\sO_Y(1)$ be a very ample invertible sheaf on $Y$ and write $\sO(n)$ for $p_2^*\sO_Y(n)$ and $\sI_Z(n)$ for $\sI_Z\otimes_{U\times_kY}\sO(n)$. Choose a point $z_i$  in each irreducible component of $p_1(z)\times_kY\setminus Z$, which is a dense open subset of $p_1(z)\times_kY$ by our assumption of the dimension of $Z$ over $U$. Since $U$ is affine, the sheaf $\sI_Z(n)$ is globally generated for all $n>>0$, so for $n>>0$,  we can find a global section $F\in H^0(Y\times_kU, \sI_Z(n))$ such that $F(z_i)\neq0$ for each $i$, and let $W_1\subset U\times_kY$ be the closed subscheme defined by $F$. By construction, $W_1$ contains $Z$ and the fiber $p_1(z)\times_UW_1$ is a codimenison one closed subscheme of $p_1(z)\times_kY$. By Chevalley's upper semi-continuity theorem \cite[Corollaire 13.1.5]{EGAIV}, after shrinking $U$ if necessary, the map $p_1:W_1\to U$ is equi-dimensional. Note that as above, the fact that $d\le \dim_kY-2$ and that $W_1$ is equi-dimenisonal of dimension $\dim_kY-1$ over $U$ implies that the fiber  $p_1(z)\times_U(W_1\setminus Z)$ is open and dense in the fiber $p_1(z)\times_UW_1$. We apply the same argument as above, replacing $U\times_kY$ with $W_1$, and find a section $G\in  H^0(Y\times_kU, \sI_Z(n'))$ for $n'>>n$ such that the subscheme $W_2$ of $U\times_kY$ has pure codimension two and is equi-dimensional  of dimension $\dim_kY-2$ over $U$. 

Let $p:W_2\to U$ be the restriction of $p_1$. Since $U\times_kY$ is smooth over $k$, $U\times_kY$ is Cohen-Macaulay, and thus $W_2$ is also Cohen-Macaulay (see \cite[Chap. 6, (16.A)(2)]{MatsumuraComAlg}).  Since $U$ is smooth over $k$, it follows that, for $u\in U$ the maximal ideal $\mathfrak{m}_u$ is generated by a regular sequence, say $x_1,\ldots, x_n$, and the fact that $W_2$ is equi-dimensional over $U$ and is Cohen-Macaulay implies that $x_1,\ldots, x_n$ defines a regular sequence in $\mathfrak{m}_w\subset \sO_{W_2,w}$ for each $w\in p^{-1}(u)$, by \cite[Chap. 6, Theorem 29]{MatsumuraComAlg}. This implies by \cite[Theorem 16.5]{MatsumuraComRingThy}  that $W_2$ is flat over $U$.  Taking a trivialization of $O(1)$ in a neighborhood of $z$ and restricting $F,G$ to $f,g\in \sI_{Z,z}$, $f,g$ is a regular sequence in $\sO_{U\times_kY,z}$ and the closed subscheme $W\subset \Spec\sO_{U\times_kY,z}$ defined by $f,g$ is flat over $\sO_{U,p_1(z)}$. 

We now discuss the case of a general $U$ of finite type over $k$; we may assume that $U$ is affine, $U=\Spec A$, for some finitely generated $k$-algebra $A$. There is thus a surjection $k[x_1,\ldots, x_n]\to A$ from some polynomial algebra, giving a closed embedding $i:U\to \A^n_k$. Let $\sI'_Z\subset \sO_{\A^n\times_kY}$ be the ideal sheaf of $(i\times\id_Y)(Z)$ and take $z\in Z$. By the case of a smooth $U$ discussed above, there is a regular sequence $F, G\in \sI'_{Z,z}$ and an open neighborhood $U'$ of $p_1(z)$ in $\A^n_k$ such that the closed subscheme $W'$ of $U'\times_kY$ defined by $F,G$ has pure codimension two in $U'\times_kY$ and is flat over $U'$. In particular, the Koszul complex
\[
0\to \sO_{U'\times_kY, z}\xrightarrow{(F,G)} \sO_{U'\times_kY, z}^2\xrightarrow{Gp_1-Fp_2}\to 
\sO_{U'\times_kY, z}\to \sO_{W',z}\to0
\]
is exact and $\sO_{W',z}$ is flat over $\sO_{U', p_1(z)}$. Thus, applying $-\otimes_{\sO_{U', p_1(z)}}\sO_{U,z}$ to this complex produces an exact sequence
\[
0\to \sO_{U\times_kY, z}\xrightarrow{(f,g)} \sO_{U\times_kY, z}^2\xrightarrow{gp_1-fp_2}\to 
\sO_{U\times_kY, z}\to \sO_{W,z}\to0
\]
where $f,g$ are the respective images of $F, G$ in $\sI_{Z,z}$ and $W=W'\times_{U'}U$. In other words, the Koszul complex for $f,g\in \sI_{Z,z}\subset \mathfrak{m}_z$ is exact. By \cite[Theorem 16.5]{MatsumuraComRingThy}, this implies that $f,g$ is a regular sequence in $\sI_{Z,z}$, which completes the proof.
\end{proof}

\begin{proposition}\label{PicProp} Let $Y$, $U$ and $Z$, $j:U\times_kY\setminus Z\to U\times_kY$  be as in Lemma~\ref{lem3}. Then the restriction map $j^*:\Pic(U\times_kY)\to \Pic(U\times_kY\setminus Z)$ is injective. 
\end{proposition}

\begin{proof} We begin by showing that, for $\sL_1, \sL_2$ invertible sheaves on $U\times_kY$, each map $s:j^*\sL_1\to j^*\sL_2$ of invertible sheaves on $U\times_kY\setminus Z$ extends uniquely to a map of invertible sheaves $\tilde{s}:\sL_1\to \sL_2$ on $U\times_kY$. To see this, note that by adjunction, $s:j^*\sL_1\to j^*\sL_2$ is equivalent to a map of $\sO_{U\times_kY}$ modules $s':\sL_1\to j_*j^*\sL_2$, and it suffices to show that $s'$ factors uniquely through the canonical map $\sL_2\to  j_*j^*\sL_2$. But for a point $z\in Z$,  choosing a local trivialization of $\sL_2$ in a neighborhood $V$ of $z$, $\sO_V\cong \sL_{2|V}$,  we see that Lemma~\ref{lem3} implies that the canonical map  $\sL_2\to  j_*j^*\sL_2$ is an isomorphism, which gives us the desired factorization. 

Now suppose we have an invertible sheaf $\sL$ on $U\times_kY$ such that $j^*\sL$ is isomorphic to $\sO_{U\times_kY\setminus Z}$ and let $s:\sO_{U\times_kY\setminus Z}\to j^*\sL$ be an isomorphism, with inverse $t:j^*\sL\to \sO_{U\times_kY\setminus Z}$. Applying the above remarks to $s,t, st$ and $ts$, we see that $s$ extends (uniquely) to an isomorphism $\tilde{s}:\sO_{U\times_kY}\to \sL$. In other words, $[j^*\sL]=0$ in $\Pic(U\times_kY\setminus Z)$ implies $[\sL]=0$ in $\Pic(U\times_kY)$, as desired. 
\end{proof}

\section{Determinants and the norm map} \label{sec:Det}
We recall some results of Deligne on the determinant of the pushforward of a locally free sheaf with respect to a finite flat morphism.

Let $f:Y\to X$ be a finite flat morphism of noetherian schemes. One has the multiplicative norm homomorphism
\[
\Nm_f:f_*\sO_Y^\times\to \sO_X^\times
\]
given by considering for $j:U\hookrightarrow X$ open and $u\in \sO_Y(f^{-1})^\times$ a unit, the determinant of the $\sO_U$-linear automorphism of locally free coherent sheaves
\[
\times u: j^*(f_*\sO_Y)\to  j^*(f_*\sO_Y)
\]

The fact that each invertible sheaf on $Y$ admits a trivializing cover of the form $f^{-1}(\sU)$ for $\sU$ a Zariski open cover of $X$ implies that this norm homomorphism induces the norm functor on the Picard groupoids
\[
\Nm_f:\sPic(Y)\to \sPic(X),
\]
compatible with the symmetric monoidal structure.

For a locally free sheaf $\sE$ on $Y$, we have the invertible sheaf $\det\sE$, and similarly for $X$. One extends the operation $\sF\mapsto \det\sF$ defined on locally free sheaves on $X$ to a homomorphism from the $K$-theory groupoid of perfect complexes $\sK(X)$ to the Picard groupoid $\sPic(X)$ 
\[
\det:\sK(X)\to \sPic(X)
\]

\begin{lemma}[\hbox{Deligne \cite[\S 7.1, (7.1.1)]{Deligne}}]\label{lem:NmAndDet} Let $f:Y\to X$ be a finite flat morphism of noetherian schemes and let $\sF$, $\sG$ be locally free coherent sheaves on $Y$ of the same rank. Then $f_*\sF$, $f_*\sG$ are locally free coherent sheaves on $X$ and there is a canonical isomorphism
\[
\det f_*\sF\cong \det f_*\sG\otimes\Nm_f(\det \sF\otimes \det\sG^{\otimes -1})
\]
\end{lemma}

\section{Spin threefolds and orientations on the Hilbert scheme of points}\label{sec:Orient}

\begin{definition} 1. Let $Y$ be a noetherian scheme. We let $D^\perf(Y)$ denote the derived category of perfect complexes on $Y$.\\[2pt]
2. Let $f:Y\to X$ be a proper flat morphism of noetherian schemes. For $\sE, \sF\in D^\perf(Y)$, we set
\[
\chi(\sE,  \sF):=Rf_*R\sHom(\sE, \sF)
\]
\end{definition}

For a  projective scheme $X$ and integer $n>0$, we have the Hilbert scheme $\Hilb^n(X)$ parametrizing flat families of closed dimension zero subschemes of length $n$ on $X$. We write this as $\Hilb^n$ when $X$ is understood.

\begin{proposition}\label{prop:main} Let $X$ be a smooth projective threefold over $k$, $n\ge1$ an integer, $\sI$ the ideal sheaf of the universal subscheme $i:Z\hookrightarrow \Hilb^n\times X$ of $\Hilb^n\times X$. Let $p_Z:Z\to \Hilb^n$ be the projection. Then 
\[
\rnk(\chi(\sI, \sI)-\chi(\sO_{\Hilb^n\times X}, \sO_{\Hilb^n\times X}))=0
\]
and
\[
\det(\chi(\sI, \sI)-\chi(\sO_{\Hilb^n\times X}, \sO_{\Hilb^n\times X}))=\det(p_{Z*}(\sO_Z))^{\otimes-2}\otimes \Nm_{Z/\Hilb^n}(i^*p_2^*\omega_{X/k})^{\otimes -1}
\]
\end{proposition}

\begin{proof}
Write $Y:=\Hilb^n\times X$ and let $i:Z\to Y$ be the inclusion. From the exact sequence
\[
0\to \sI\to \sO_Y\to i_*\sO_Z\to 0
\]
we get (in $\sK(\Hilb^n)$)
\[
\chi(\sI,\sI)=\chi(\sO_Y,\sO_Y)-\chi(\sO_Y,i_*\sO_Z)-
\chi(i_*\sO_Z,\sO_Y)+\chi(i_*\sO_Z,i_*\sO_Z).
\]

We use Grothendieck-Serre duality to compute $\chi(i_*\sO_Z,\sO_Y)$.   Since $p_Z$ and $i$ are finite, we have  
\[
p_{Z*}=Rp_{Z*}=Rp_{1*}\circ Ri_*=Rp_{1*}\circ i_*
\]
and  we have
\begin{align*}
Rp_{1*}\RsHom_{\sO_Y}(i_*\sO_Z, \sO_Y)&= Rp_{1*}\RsHom_{\sO_Y}(i_*\sO_Z\otimes p_2^*\omega_{X/k}, p_2^*\omega_{X/k})\\
&= Rp_{1*}\RsHom_{\sO_Y}(i_*(\sO_Z\otimes i^*p_2^*\omega_{X/k}), p_1^!\sO_{\Hilb^n}[-3])\\
&=R\sHom_{\sO_{\Hilb^n}}(Rp_{1*}\circ i_*(\sO_Z\otimes i^*p_2^*\omega_{X/k}), \sO_{\Hilb^n})[-3]\\
&=R\sHom_{\sO_{\Hilb^n}}(p_{Z*}(\sO_Z\otimes i^*p_2^*\omega_{X/k}), \sO_{\Hilb^n})[-3].
\end{align*}
Thus
\[
\det(Rp_{1*}\RsHom_{\sO_Y}(i_*\sO_Z, \sO_Y))=\det(p_{Z*}(\sO_Z\otimes i^*p_2^*\omega_{X/k})).
\]

Since $Z\to \Hilb^n$ is flat and finite of degree $n$, $p_{Z*}(\sO_Z\otimes i^*p_2^*\omega_{X/k})$ is locally free of rank $n$, so 
\[
\rnk Rp_{1*}\RsHom_{\sO_Y}(i_*\sO_Z, \sO_Y)=(-1)^3n= -n.
\]

Similarly
\begin{align*}
\det\chi(\sO_Y, i_*\sO_Z)&=\det Rp_{1*}\RsHom(\sO_Y, i_*\sO_Z)\\
&=\det Rp_{1*}(i_*\sO_Z)\\
&=\det p_{Z*}(\sO_Z),
\end{align*}
and 
\[
\rnk Rp_{1*}(i_*\sO_Z)=n.
\]

By Lemma~\ref{lem:NmAndDet}, we have 
\[
\det(p_{Z*}(\sO_Z\otimes i^*p_2^*\omega_{X/k}))=\det(p_{Z*}(\sO_Z))\otimes \Nm_{Z/\Hilb}( i^*p_2^*\omega_{X/k}).
\]

For the term $\chi(i_*\sO_Z, i_*\sO_Z)$, since $Y\to \Hilb^n$ is smooth and $p_Z$ is flat, $i_*\sO_Z$ admits a finite resolution by locally free sheaves on $Y$,  $\sE_*\to i_*\sO_Z\to 0$, and 
\[
Rp_{1*}\RsHom_{\sO_Y}(i_*\sO_Z, i_*\sO_Z)=Rp_{1*}\sHom_{\sO_Y}(\sE_*,i_*\sO_Z)=
p_{Z*}(i^*\sE_*^\vee)
\]
Let $j:Y\setminus Z\to Y$ be the inclusion. Then  $j^*\sE^\vee_*$ is exact, thus 
$\det(j^*\sE^\vee_*)=\sO_{Y\setminus Z}$. By Proposition~\ref{PicProp}, this implies that $\det(\sE^\vee_*)=\sO_Y$. Clearly $\sE_*^\vee$ has zero rank, so $i^*\sE_*^\vee$ has rank zero and determinant $\sO_Z$. By Lemma 2.3.8, we thus  have
\begin{align*}
\det\chi(i_*\sO_Z, i_*\sO_Z)&=\det Rp_{1*}\RsHom_{\sO_Y}(i_*\sO_Z, i_*\sO_Z)\\
&=\det p_{Z*}(i^*\sE_*^\vee)\\
&=\Nm_{Y/\Hilb^n}(\det i^*\sE_*^\vee)\\
&=\sO_{\Hilb^n}
\end{align*}
and
\[
\rnk\chi(i_*\sO_Z, i_*\sO_Z)=0
\]

Putting this all together gives 
\[
\rnk(\chi(\sI, \sI)-\chi(\sO_Y, \sO_Y))=0
\]
and
\[
\det(\chi(\sI, \sI)-\chi(\sO_Y, \sO_Y))=\det(p_{Z*}(\sO_Z))^{\otimes -2}\otimes\Nm(i^*p_2^*\omega_{X/k})^{\otimes -1}
\]
\end{proof}

\section{The main theorem}\label{sec:Main}

\begin{definition} Let $X$ be a smooth proper scheme over a field $k$. A {\em spin structure} on $X$ is a choice of an invertible sheaf $\sL$ on $X$ and an isomorphism $\rho:\omega_{X/k}\xrightarrow{\sim} \sL^{\otimes 2}$.
\end{definition}

Let $X$ be a smooth projective threefold and fix an integer $n>0$, giving us the Hilbert scheme $\Hilb^n$. Recall from \cite[Lemma 1.2.2]{BehrendFantechi}\cite{Thomas} the perfect obstruction theory $\phi_\bullet:E_\bullet\to L_{\Hilb^n(X)/k}$.  The complex $E_\bullet$  is by definition
\begin{equation}\label{eqn:POT}
E_\bullet= Rp_{1*}(\sF\otimes p_2^*\omega_{X/k})[2]\cong Rp_{1*}(\RsHom(\sF, p_2^*\omega_{X/k})[2]),
\end{equation}
with $\sF$ fitting into a distinguished triangle
\begin{equation}\label{eqn:POTTri}
\sF\to \RsHom(\sI_Z, \sI_Z)\xrightarrow{\Tr} \sO_{Y}\to \sF[1],
\end{equation}
where $\Tr$ is the trace map. Following \cite{HT}, the map $\phi_\bullet:E_\bullet\to L_{\Hilb^n(X)/k}$ arises from the relative
 Atiyah class $At_\sI: \RsHom(\sI_Z, \sI_Z)\to p_2^*\Omega_{X/k}[1]$ (see details in \S\ref{sec:equi}).

\begin{theorem} \label{thm:Main} Let $X$ be a smooth projective threefold over a field $k$, endowed with a spin structure $\rho:\omega_{X/k}\xrightarrow{\sim} \sL^{\otimes 2}$.   Let $i:Z\to \Hilb^n(X)\times_kX$ be the universal subscheme and let $E_\bullet\in D^\perf(Z)$ be as in \eqref{eqn:POT}. Then  $E_\bullet$ has virtual rank zero and admits a canonical isomorphism 
\[
\rho_{\Hilb}: \det E_\bullet\to [\det(p_{Z*}(\sO_Z))\otimes\Nm_{Z/\Hilb^n(X)}(i^*p_2^*\sL)]^{\otimes -2}
\]
\end{theorem}

\begin{proof} As before, write $Y:= \Hilb^n\times X$. 
By Serre duality, the distinguished triangle \eqref{eqn:POTTri} gives the distinguished triangle
\begin{multline*}
E_\bullet\to \RsHom_{\sO_{\Hilb^n}}(Rp_{1*}\RsHom(\sI_Z, \sI_Z), \sO_{\Hilb^n})[-1]\\\to
 \RsHom_{\sO_{\Hilb^n}}(Rp_{1*}\sO_{Y}, \sO_{\Hilb^n})[-1]\to E_\bullet[1].
 \end{multline*}
This in turn gives us a canonical isomorphism in $\sK(X)$
\[
[E_\bullet]\cong [\RsHom_{\sO_{\Hilb^n}}(\chi(\sI_Z, \sI_Z)-\chi(\sO_{Y}, \sO_{Y}), \sO_{\Hilb^n})][-1],
\]
 giving the isomorphism  in $\Pic(X)$
\[
\det E_\bullet\cong \det(\chi(\sI_Z, \sI_Z)-\chi(\sO_{Y}, \sO_{Y})).
\]
and the identity
\[
\rnk E_\bullet =\rnk(\chi(\sI_Z, \sI_Z)-\chi(\sO_{Y}, \sO_{Y})).
\]

By Proposition~\ref{prop:main}, we have
\[
\rnk E_\bullet =0
\]
and we have a canonical isomorphism
\[
\det E_\bullet\cong \det(p_{Z*}(\sO_Z)^{\otimes -2}\otimes\Nm_{Z/\Hilb^n(X)}(i^*p_2^*\omega_{X/k})^{\otimes -1};
\]
inserting the isomorphism $\rho:\omega_{X/k}\xrightarrow{\sim} \sL^{\otimes 2}$ gives the result.
\end{proof}

\section{The equivariant case}\label{sec:equi} Let $G$ be an affine  group scheme, smooth and of finite type over $k$.  We suppose that $X$ is endowed with a $G$-action; this induces a $G$-action on $\Hilb^n(X)$, with $Z\subset \Hilb^n(X)\times X$ a $G$-stable closed subscheme.

Recall that a noetherian scheme $W$ with an action of an algebraic group scheme $\sG$ has the {\em $\sG$-resolution property} if each $\sG$-linearized coherent sheaf $\sF$ on $W$ admits a $\sG$-equivariant surjection $\sV\to \sF$, with $\sV$ a $\sG$-linearized locally free coherent sheaf on $W$.

\begin{lemma}\label{lem:GLinearization0}  Suppose that there is a homomorphism $G\to H$, with $H$ a connected linear algebraic group over $k$, such that the $G$-action on $X$ extends to an $H$-action. Then $X$ admits a $G$-linearized very ample invertible sheaf.
\end{lemma}

\begin{proof} By  \cite[Chap. 1, \S3, Proposition 1.5, Corollary 1.6]{GIT}, $X$ has a very ample invertible sheaf with an $H$-linearization; the result follows immediately.
\end{proof}

\begin{lemma}\label{lem:GLinearization} Suppose there is a  very ample $G$-linearized invertible sheaf on $X$. Then \\[5pt]
1.  $\Hilb^n(X)$ and $\Hilb^n(X)\times X$ admit respective $G$-linearized very ample invertible sheaves. \\[2pt]
2. $X$, $ \Hilb^n(X)$ and $\Hilb^n(X)\times X$ all have the $G$-resolution property.
\end{lemma}

\begin{proof} Applying \cite[Lemma 2.6.]{Thomason}, we see that (2) follows from (1).

For (1), letting $V=H^0(X, \sO_X(1))$, we have the canonical closed embedding $i:X\to \P(V)$, $\P(V)$ inherits a $G$-action from the induced $G$-action on $V$ and $i$ is $G$-equivariant. In addition, $\sO_{\P(V)}(1)$ has a canonical $G$-linearization. By its construction, the Hilbert scheme $\Hilb^n(X)$ is  a closed subscheme of a Grassmannian $\Gr(N, \Sym^M(V)^\vee)$ for certain integers $N,M>0$, and the action of $G$ on $\Hilb^n(X)$, induced naturally from the $G$-action on $X$ via its representing the functor of length $n$ dimension zero closed subschemes of $X$, is the restriction of the $G$-action on $\Gr(N, \Sym^M(V)^\vee)$ induced from the $G$-action on $V$. Restricting the Pl\"ucker embedding $\Gr(N, \Sym^M(V)^\vee)\hookrightarrow \P(\Lambda^N\Sym^M(V)^\vee)$ to $\Hilb^n(X)$, we see that the restriction of $\sO_{\P(\Lambda^N\Sym^M(V)^\vee)}(1)$ to $\Hilb^n(X)$ is a $G$-linearized very ample invertible sheaf on $\Hilb^n(X)$. Taking the external tensor product of the respective $G$-linearized invertible sheaves on $X$ and on $\Hilb^n(X)$ gives us a $G$-linearized very ample invertible sheaf on $X\times\Hilb^n(X)$.  This completes the proof of (1).
\end{proof}

For the remainder of this section, we assume that there is a very ample $G$-linearized invertible sheaf on $X$.

We note that the perfect complex $E_\bullet$ admits a canonical $G$-linearization. To see this, we choose a $G$-equivariant resolution of $\sI$ by locally free $G$-linearized sheaves $\sP_\bullet\to \sI$, which exists by Lemma~\ref{lem:GLinearization}. In fact there is a finite such resolution, since $\sI$ is flat over $\Hilb^n(X)$ and $X$ is smooth over $k$. The trace map $\RsHom(\sI, \sI)\to \sO_Y$ is represented by the map of complexes $\Tr_\sP: \sHom(\sP_\bullet, \sP_\bullet)\to \sO$ which is just the sum of the trace maps $\sHom(\sP_i,\sP_i)\to \sO$. We can thus take $\sF:=\Cone(\Tr)[-1]$; since $\Tr_\sP$ is $G$-equivariant, this gives $\sF$ a $G$-linearization. Once we have this, $\sHom(\sF, \omega_{Y/\Hilb^n})[2]$ inherits a canonical $G$-linearization from $\sF$ and the canonical $G$-linearization on $\omega_{Y/\Hilb^n}$. Since $\sF$ is a complex of locally free coherent sheaves on $Y$, we have
\[
E_\bullet:=Rp_{1*}\RsHom(\sF, \omega_{Y/\Hilb^n})[2]=Rp_{1*}\sHom(\sF, \omega_{Y/\Hilb^n})[2]
\]
hence $E_\bullet$ inherits a $G$-linearization from the one we have just constructed on $\sHom(\sF, \omega_{Y/\Hilb^n})[2]$.   

The isomorphism in Lemma~\ref{lem:NmAndDet} is natural in $f$, $\sF$ and $\sG$, hence extends immediately to a $G$-equivariant isomorphism,  and the argument used in Proposition~\ref{prop:main} goes through without change to yield the equivariant version. For Theorem~\ref{thm:Main}, one needs a ``$G$-linearized spin structure'' on $X$, that is, a $G$-linearized invertible sheaf with a $G$-equivariant isomorphism
$\rho:\omega_{X/k}\xrightarrow{\sim} \sL^{\otimes 2}$, where one gives $\omega_{X/k}$ the natural $G$-linearization induced from the $G$-action on $X$. For the sake of completeness, we state the equivariant version of Theorem~\ref{thm:Main} here. 

\begin{theorem} \label{thm:MainEquiv} Let $G$ be a finite type group-scheme over $k$. Let $X$ be a smooth projective threefold over a field $k$, endowed with $G$-action. Suppose there is a $G$-linearized very ample invertible sheaf $\sL$ on $X$ and we are given a $G$-linearized spin structure $\rho:\omega_{X/k}\xrightarrow{\sim} \sL^{\otimes 2}$. Let $i:Z\to \Hilb^n(X)\times_kX$ be the universal subscheme.  We endow $\RsHom(\sI_Z, \sI_Z)$ with a $G$-linearization as above. Let  $E_\bullet\in D^\perf(Z)$ be as given by \cite[Lemma 1.2.2]{BehrendFantechi}, with its canonical $G$-linearization induced from the given one on $\RsHom(\sI_Z, \sI_Z)$, as described above. Then the complex $E_\bullet$ has virtual rank zero and there is a canonical $G$-equivariant isomorphism 
\[
\rho_{\Hilb}: \det E_\bullet\to [\det(p_{Z*}(\sO_Z))\otimes\Nm_{Z/\Hilb^n(X)}(i^*p_2^*\sL)]^{\otimes -2}
\]
\end{theorem}

We conclude with a few words about enhancing the Thomas perfect obstruction theory on $\Hilb^n(X)$ to a $G$-equivariant version. This is already discussed by Behrend-Fantechi in \cite[\S 2]{BehrendFantechi}. The $G$-linearization of $E_\bullet$ is constructed there, essentially as described above, but there is no explicit description of the map $\phi:E_\bullet\to L_{\Hilb^n(X)/k}$, other than to refer to the work of Thomas \cite{Thomas}. However in that work, the virtual fundamental class is constructed using the obstruction sheaf-tangent sheaf furnished by the cohomology of $E_\bullet^\vee$. There seems to be no mention of the cotangent complex in \cite{Thomas}; probably it is clear to the experts how to pass from the obstruction sheaf-tangent sheaf data to a construction of the map $\phi$. Instead, 
we refer to  subsequent work of Huybrechts-Thomas \cite{HT} that gives an explicit construction of $\phi$ using their truncated Atiyah class. The equivariant case is discussed in detail by Ricolfi \cite{Ricolfi}, building on the work in \cite{HT}.

In our case, since we are working with the smooth morphism $p_1:\Hilb^n(X)\times_kX\to X$, we can use the usual universal Atiyah class for $X$, pulled back to $\Hilb^n(X)\times_kX$, to construct the map $At_Z: \RsHom(\sI_Z,\sI_Z)\to  L_{\Hilb^n(X)\times X/\Hilb^n(X)}[1]$ as follows.  Let $i:X\to X\times_kX$ be the diagonal and $\sI_{\Delta_X}$ the ideal sheaf of the diagonal. The universal Atiyah class for $X$ is the map $At_X:i_*\sO_X\to i_*\Omega_{X/k}[1]$ in $D^\perf(X\times_kX)$ given by  the extension of $\sO_{X\times_kX}$-modules
\[
0\to\sI_X/\sI_X^2\to \sO_{X\times_kX}/\sI_X^2\to i_*\sO_X\to0
\]
together with the canonical isomorphism $i_*\Omega_{X/k}\cong \sI_X/\sI_X^2$. Let $p_1, p_2:\Hilb^n(X)\times X\times X\to \Hilb^n(X)\times_kX$, $p_X:\Hilb^n(X)\times X\to X$ and 
$p_{X\times X}:\Hilb^n(X)\times X\times \to X\times X$ be the projections. For a perfect complex $\sP$ on $\Hilb^n(X)\times_kX$, we have the induced map
\[
At_\sP:\sP\to \sP\otimes_k L_{\Hilb^n(X)\times_kX/\Hilb^n(X)}[1]
\]
defined by  $Rp_{2*}(p_1^*(\id_\sP\otimes p_{X\times X}^*At_X)$, that is, we apply $Rp_{2*}$ to the map
\[
p_1^*\sP\otimes p_{X\times X}^*i_*\sO_X\xrightarrow{\id_\sP\otimes p_{X\times X}^*At_X} (p_1^*\sP\otimes p_{X\times X}^*i_*\Omega_{X/k}
\]
and use the identifications 
\[
p_X^*\Omega_{X/k}\cong p_X^*L_{X/k}\cong L_{\Hilb^n(X)\times_kX/\Hilb^n(X)}
\]
\[
Rp_{2*}(p_1^*\sP\otimes p_{X\times X}^*i_*\sO_X)\cong \sP
\]

We let 
\[
At'_\sP:\RsHom(\sP, \sP)=\sP^\vee\otimes\sP\to L_{\Hilb^n(X)\times_kX/\Hilb^n(X)}[1]
\]
be the map induced by $At_\sP$ by the adjunction $\Hom(\sP^\vee\otimes A, B)\cong \Hom(A, \sP\otimes B)$.

Following \cite{HT}, we have  map $\phi:E_\bullet\to L_{\Hilb^n(X)/k}$ in $D^-(\Hilb^n(X))$ 
given by the composition
\begin{multline*}
Rp_{1*}(\sF\otimes p_2^*\omega_{X/k}[2])\to
Rp_{1*}(\RsHom(\sI_Z, \sI_Z)\otimes p_2^*\omega_{X/k}[2])\\\xrightarrow{Rp_{1*}(At'_{\sI_Z}\otimes\id)} Rp_{1*}(L_{\Hilb^n(X)\times X/\Hilb^n(X)}\otimes p_1^!\sO_{\Hilb^n(X)})\\
\xrightarrow{\del} Rp_{1*}(p_1^*L_{\Hilb^n(X)/k}\otimes p_1^!\sO_{\Hilb^n(X)})\xrightarrow{\Tr}
L_{\Hilb^n(X)/k},
\end{multline*}
Here we use the Atiyah map
\[
\del:L_{\Hilb^n(X)\times X/\Hilb^n(X)}\to L_{\Hilb^n(X)/k}[1]
\]
and the trace map 
\begin{multline*}
Rp_{1*}(p_1^*L_{\Hilb^n(X)/k}\otimes p_1^!\sO_{\Hilb^n(X)})\\\cong
L_{\Hilb^n(X)/k}\otimes Rp_{1*}(p_1^!\sO_{\Hilb^n(X)})   \xrightarrow{\id\otimes\Tr}
L_{\Hilb^n(X)/k}
\end{multline*}
 By \cite[Theorem 4.1, Corollary 4.3]{HT}, this defines a perfect obstruction theory; it follows from the proof of {\it loc. cit.} that this recovers the perfect obstruction theory on $\Hilb^n(X)$ constructed by Thomas \cite{Thomas}.

 If one works over $\C$, it is shown in \cite{Ricolfi} that this admits a canonical extension to the $G$-equivariant setting. With minor modifications to ensure of smoothness of the given group scheme $G$ over $k$, and replacing the  assumption that $G$ is reductive with linear reductivity, 
 the arguments of {\it loc. cit.}, together with  Lemma~\ref{lem:GLinearization} show the following.

\begin{theorem} Let $k$ be a field and let $G$ be a linearly reductive smooth affine algebraic group scheme over $k$. Let $X$ be a smooth projective 3-fold over $k$, endowed with a $G$-action and a $G$-linearized very ample invertible sheaf, and let $D^-_G(\Hilb^n(X))$ denotes the (bounded above) derived category of $G$-linearized coherent sheaves on $\Hilb^n(X)$. Let $E_\bullet:=Rp_{1*}(\sF\otimes p_2^*\omega_{X/k}[2])$, endowed with the $G$-linearization described above and give $L_{\Hilb^n(X)}$ its canonical $G$-linearization. Then there is a  morphism $\phi^G:E_\bullet\to L_{\Hilb^n(X)}$ in $D^-(\Hilb^n(X))$ that defines a $G$-equivariant perfect obstruction theory and yields the Thomas perfect obstruction theory by forgetting the $G$-action.
\end{theorem}

\begin{proof} The Atiyah map $AT_\sP$ extends to give a morphism 
\[
AT^G_\sP: \sP\to \sP\otimes_k L_{\Hilb^n(X)\times_kX/\Hilb^n(X)}[1]
\]
in $D^\perf_G(\Hilb^n(X)\times_kX)$ for $\sP\in D^\perf_G(\Hilb^n(X)\times_kX)$ (natural in $\sP$).   Using Lemma~\ref{lem:GLinearization} realizes $\sI_Z$ as an object in $D^\perf_G(\Hilb^n(X)\times X)$, which yields the $G$-equivariant lifting $\phi^G$ of $\phi$.
\end{proof}

\end{document}